\documentclass[a4paper,11pt]{amsart}
\usepackage[latin1]{inputenc}
\usepackage[english]{babel}
\usepackage{amsmath}
\usepackage{amsfonts}
\usepackage{amsthm}
\usepackage{amssymb}

\newcommand{\Z}{\mathbb{Z}}

\newcommand{\C}{\mathbb{C}}

\newcommand{\N}{\mathbb{N}}

\newcommand{\del}{\partial}
\newcommand{\delbar}{\overline\partial}
\newcommand{\norm}[1]{\left\lVert#1\right\rVert}
\newcommand{\Cinf}{\mathcal{C}^{\infty}}
\newcommand{\Ha}{\mathcal{H}}

\DeclareMathOperator{\Ker}{Ker}
\DeclareMathOperator{\Ima}{Im}
\DeclareMathOperator{\Dim}{Dim}

\newtheorem{teo}{Theorem}

\newtheorem{prop}{Proposition}

\newtheorem{oss}{Remark}
\newtheorem{lemma}{Lemma}
\newtheorem{hyp}{Assumption}

\begin{document}
\title{On the degeneration of the Fr\"olicher spectral sequence and small deformations}
\author{Michele Maschio}

\begin{abstract}
We study the behavior of the degeneration at the second step of the Fr\"olicher spectral sequence of a $\Cinf$ family of compact complex manifolds. Using techniques from deformation theory and adapting them to pseudo-differential operators we prove a result \textit{\`a la Kodaira-Spencer} for the dimension of the second step of the Fr\"olicher spectral sequence and we prove that, under a certain hypothesis, the degeneration at the second step is an open property under small deformations of the complex structure. We also provide an example in which we show that, if we omit that hypothesis, we lose the openness of the degeneration.
\end{abstract}

\maketitle

\begin{section}{Introduction}
An important topic in complex geometry is the study of the behavior of certain properties under $\Cinf$ changes of the complex structure of a compact complex manifold. In one of the most celebrated papers, \cite{kodaira1960deformations},the stability of the K\"ahler condition has been proved, namely if a compact complex manifold admits a K\"ahler metric and we change in $\Cinf$ way its complex structure, then on the new complex manifold can be defined a K\"ahler metric.
In this paper we study the stability under small deformations of the complex structure of the degeneration at the second step of the Fr\"olicher spectral sequence. Spectral sequences are generalizations of exact sequences and they were introduced by Leray in \cite{leray1950lann} to compute the homology groups of sheaves by taking successive approximations.  Given a compact complex manifold $(M,J)$, the Fr\"olicher spectral sequence $(E_r^{\bullet,\bullet},d_r)$ is the the spectral sequence of the double complex $(\Lambda^{p,q},\del,\delbar)$ of $\Cinf$ $(p,q)$-forms of $(M,J)$ (see \cite{frolicher1955relations}). It provides a link between the Dolbeault and the de Rham cohomology groups. Moreover, it "measures" the failure of results in cohomology theory that are valid for K\"ahler manifolds (or more generally for manifolds satisfying the $\del\delbar$-Lemma, see \cite{deligne1975real} ). 
The degeneration at the second step of the Fr\"olicher spectral sequence has been studied recently in \cite{popovici2016degeneration}; the author constructed a pseudo differential operator $\tilde\Delta$ whose kernel is isomorphic to $E_2^{\bullet,\bullet}$.
Since this operator also provides relations between SKT and Gauduchon metrics and super SKT and strongly Gauduchon metrics  of $(M,J)$ (see \cite{gauduchon1977theoreme}, \cite{popovici2009limits}, \cite{gauduchon19841} and \cite{maschio2017cone} for more information about those metrics) and since, always in \cite{popovici2016degeneration}, Popovici conjectured a connection between the degeneration at the second step of $(E^{\bullet,\bullet}_r,d_r)$ and the presence of an SKT metric on $(M,J)$,  in the second section we prove some properties of $\tilde\Delta$.

As stated above, we are interested in the study the behavior of the degeneration at the second step of $(E^{\bullet,\bullet}_r,d_r)$ on a compact complex manifold $(M,J)$ under small deformations $J_t$ of the complex structure $J=J_0$. More precisely, under the assumption
\begin{hyp}\label{ass}
For every $(p,q)\in \Z^2$, the Hodge numbers $h^{p,q}(t)$, i.e. the dimension of the Dolbeault cohomology groups of $(M,J_t)$,  are constants with respect to $t$.
\end{hyp}
where $\Delta''_t$ is the Dolbeault Laplacian of $(M,J_t)$, we prove the following
\begin{teo}\label{MT}
Let $(M,J_t)$ be a family of compact complex manifolds. If, for every $(p,q)\in\Z^2$, $h^{p,q}(t)$ is independent of $t$ and $E^{p,q}_2(0)\simeq E^{p,q}_{\infty}(0)$. Then, for every $(p,q)\in\Z^2$, $E^{p,q}_2(t)\simeq E^{p,q}_{\infty}(t)$ for $t$ sufficiently closed to $0$. 
\end{teo}

The basic tools in the proof of Theorem \ref{MT} are the following: the general \textit{a priori estimate} 
\begin{teo}\label{ape}
For every $k\in\Z$ there exists a constant $C_k$ depending only on $k$ such that for every $\phi\in \Lambda^{p,q}(M)$
\begin{equation}
\norm{\phi}_{k+2}^2\leq C \left( \norm{\tilde{\Delta}\phi}^2_{k}+\norm{\phi}^2_0 \right) 
\end{equation}
\end{teo}
and a property of self-adjoint elliptic differential operators that also holds for $\tilde\Delta$, namely
\begin{teo}\label{ker}
Fix a $\Cinf$ family $\{\omega_t\}$ of Hermitian metrics. If  $\Dim \Ker \Delta''_t$ is independent of $t$, then $Dim \Ker \tilde\Delta_t$ is an upper semi-continuous function in $t$.
\end{teo}

The paper is organized as follows:

In Section \ref{pre} we fix notations and we recall basic facts on the Fr\"olicher spectral sequence, differential operators and Sobolev norms.

In Section \ref{nl} we recall the construction of $\tilde\Delta$ recently made by Popovici in \cite{popovici2016degeneration} and we give the proof of Theorem \ref{ape}.

Section \ref{defo} is devoted to the proof Theorem \ref{MT}. First of all using techniques \`a la Kodaira \cite{kodaira2006complex}, we show Theorem \ref{ker} applied to the family of self-adjoint "elliptic" pseudo-differential operators $\{\tilde\Delta_t\}$. Then we use Theorems \ref{ape} and \ref{ker} to prove Theorem \ref{MT}.

In the last section, starting with the completely solvable Nakamura manifold and taking suitable complex deformations of the complex structure studied in  \cite{tomassini2014dolbeault}, we show, by direct computations, that if Assumption \ref{ass} is not satisfied, then also the conclusion of Theorem \ref{MT} does not hold.
\end{section}

\begin{section}{Preliminaries}\label{pre}
Let $(M,J)$ be a compact complex manifold of complex dimension $n$. The main object of this paper is the \textit{Fr\"olicher Spectral Sequence} $(E^{p,q}_r,d_r)$ of $(M,J)$. Namely, denoting by $(\Lambda^{\bullet,\bullet}(M),\del,\delbar)$ the double complex of $\mathcal{C}^{\infty}$ forms over $M$, the Fr\"olicher spectral sequence $(E^{p,q}_r,d_r)$ of $(M,J)$ is the spectral sequence associated to $(\Lambda^{\bullet,\bullet}(M),\del,\delbar)$ and it is constructed in the following way:
let  $ E^{p,q}_0:=\Lambda^{p,q}(M)$ and  $d_0:=\delbar:\Lambda^{p,q}(M)\rightarrow \Lambda^{p,q+1}(M)$, then, inductively, let  $E^{p,q}_r$ be the cohomology group of the $(r-1)$-th step of the sequence and $d_r:E^{p,q}_r\rightarrow E^{p+r,q+r-1}_r$.  This spectral sequence  satisfies the following properties:
\begin{itemize}
\item $E^{p,q}_1$ is isomorphic to the $(p,q)$-th group of the Dolbeault cohomology;
\item for every $r\geq 1$, $E^{p,q}_r$ is a finite-dimensional complex vector space;
\item for every $r\geq 1$, $\dim E^{p,q}_r\geq \dim E^{p,q}_{r+1}$.
\end{itemize}
Moreover the sequence is said to \textit{degenerate at the step r} if $r$ is the smallest integer such that, for every $(p,q)\in\Z^2$ and every $r'>r$, $\dim E^{p,q}_r= \dim E^{p,q}_{r'}$, when this happens we have that $H^k_{dR}(M)\simeq \oplus_{p+q=k} E^{p,q}_r$. Since the sequence degenerates at the first step when $M$ is a K\"ahler manifold, we can say that it provides a "measure of the non-K\"ahlerianity" of a manifold. In particular we are interested in the degeneration at the second step of $(E^{p,q}_r,d_r)$.

Although the spectral sequence is an algebraic object, we use analytic arguments to prove Theorem \ref{ker}. This is a classic approach, in fact, for example, the \textit{Dolbeault Laplacian} can be used to study the Dolbeault cohomology of a compact complex manifold. Since we will use it later on, we recall its definition: fix a Hermitian metric $g$ over $M$ and let $*$ be the Hodge-star operator with respect to $g$, we recall that the Dolbeault Laplacian $\Delta''$ is defined as 
\begin{equation}\label{dl}
\Delta'':=\delbar\delbar^*+\delbar^*\delbar,
\end{equation}
where $\delbar^*:=-*\del*$ is the operator adjoint to $\delbar$. It is a standard result that $\Delta''$ is an elliptic, self-adjoint and non-negative differential operator. Moreover $\Delta''$ induces the decompositions
$$
\Lambda^{p,q}(M)=\Ker \Delta'' \oplus \Ima \Delta''=\Ker \Delta'' \oplus \Ima \delbar \oplus \Ima \delbar.
$$
For every $(p,q)$, we denote with $\mathcal{H}^{p,q}:=\Lambda^{p,q}(M)\cap \Ker \Delta''$ the space of $\Delta''$-harmonic $(p,q)$ forms. It is a well-known fact that $\Ha^{p,q}$ is a finite dimensional $\C$-vector space.

Finally, since we need some estimates, we recall the definition of the \textit{Sobolev norm} of a $(p,q)$-form. Let $\{ U_j\}$ be a finite covering by coordinate neighborhoods of $M$, let $\{x^i_j\}$ be coordinate over $U_j$ and let $\{\eta_j\}$ be a partition of unity subordinated to $\{U_j\}$. Given a $\Cinf$ function $f$ over $M$, the $k$-th Sobolev norm of $f\in \Cinf(M,\C)$ is defined as
\begin{equation}\label{SN}
\norm{f}_k^2:= \sum_{|L|=0}^k \sum_j \sum_{D_j^l} \int_{U_j} |D_j^L f_j(x)|^2 d X_j,
\end{equation}
where
\begin{itemize}
\item $f_j:=\eta_j f$;
\item $L=(l_1,\dots, l_{2n})$ is a multiindex and $|L|=l_1+\dots l_{2n}$;
\item $D_j^L:=\frac{\del^{l_1}}{\del x_j^{1,l_1}}\dots\frac{\del^{l_{2n}}}{\del x_j^{2n,l_{2n}}}$ is a differential operator of rank $|L|$;
\item $dX_j$ is the Lebesque measure associated to $g$ and expressed in the local coordinates $dx_j^{1},...,dx_j^{2n}$.
\end{itemize}

If $\phi$ is an $r$-form over $M$ that, for every $U_j$, can be written locally as $\phi=\sum_{A_l} f_j^{A_l} dx^{A_l}$ ($dx^{A_l}$ stands for $dx_{l_1}\wedge...\wedge dx_{l_r}$ with $\{l_1,...,l_k\}=A_l$ and $l_1<...<l_k$), then we define the $k$-th \textit{Sobolev norm} of $\phi$ as 
\begin{equation}\label{SNF}
\norm{\phi}_k^2:=\sum_{A_l}\norm{f^{A_l}}^2_k.
\end{equation}

\end{section}

\begin{section}{The new Laplacian $\tilde{\Delta}$}\label{nl}

In \cite{popovici2016degeneration}, Popovici introduced a pseudo-differential operator related to the second step of the Fr\"olicher spectral sequence in the same way the Dolbeault Laplacian is related to the first step, i.e., the Dolbeault cohomology of the manifold.

We recall its construction: let
\begin{equation}\label{sl}
\Delta'_{p''}:=\del p''\del^* +\del^* p''\del,
\end{equation}
where, for every $(p,q)$, $p'':\Lambda^{p,q}(M)\rightarrow \Ha^{p,q}$ is the orthogonal projection and $\del^*:=-*\delbar *$ is the adjoint operator of $\del$. 
Then Popovici defined the pseudo-differential operator
\begin{equation}\label{fl}
\tilde{\Delta}:=\Delta'_{p''}+\Delta'':\Lambda^{p,q}(M)\rightarrow \Lambda^{p,q}(M),
\end{equation}
for every $(p,q)\in \Z^2$.
$\tilde{\Delta}$ is not a differential operator, but it is still possible to prove that it satisfies properties of elliptic operator. In particular we have
\begin{teo}
For all $p,q$, $\tilde{\Delta}:\Lambda^{p,q}(M)\rightarrow \Lambda^{p,q}(M)$ behaves like an elliptic self-adjoint differential operator in the sense that $\Ker \tilde{\Delta}$ is a finite dimensional $\C$-vector space, $\Ima \tilde{\Delta}$ is closed and finite codimensional in $\Lambda^{p,q}(M)$, there is an orthogonal (for the $L^2$ inner product induced by $g$) 2-space decomposition
$$
\Lambda^{p,q}(M)=\Ker\tilde\Delta \oplus \Ima \tilde\Delta
$$
giving rise to an orthogonal 3-space decomposition
$$
\Lambda^{p,q}(M)=\Ker\tilde\Delta \oplus \left( \Ima \delbar + \Ima\del_{\mid_{\Ker \delbar}}\right) \oplus \left( \Ima  \delbar^*+ \Ima (\del^* \circ p'')\right)
$$
in which 
$$
\begin{array}{l}
\Ker\tilde\Delta \oplus \left( \Ima \delbar + \Ima \del_{\mid_{\Ker \delbar}}\right)=\Ker (p''\circ \del)\cap \Ker \delbar;\\
\Ker\tilde\Delta  \oplus \left( \Ima  \delbar^*+ \Ima (\del^* \circ p'')\right) = \Ker (p''\circ\del^*)\cap \Ker\delbar^* ;\\
\left( \Ima \delbar + \Ima\del_{\mid_{\Ker \delbar}}\right) \oplus \left( \Ima  \delbar^*+ \Ima (\del^* \circ p'')\right) = \Ima \tilde\Delta.
\end{array}
$$
Moreover, $\tilde\Delta$ has a compact resolvent which is a pseudo-differential operator $G$ of order $-2$, the Green operator of $\tilde\Delta$, hence the spectrum of $\tilde\Delta$ is discrete and consists of non-negative eigenvalues that tend to $+\infty$.
\end{teo} 

In the same paper, Popovici constructed an isomorphism between the kernel of $\tilde\Delta\cap \Lambda^{p,q}(M,J)$ and $E^{p,q}_2$, for every $(p,q)\in\Z^2$.

Given the previous theorems, it is natural to study the behavior of $\tilde\Delta$ in order to better understand the Fr\"olicher spectral sequence. In particular we want to study this behavior under small deformations of the complex structure, but we need some preliminary results.

First we study some norm estimates for $\tilde\Delta$. We observe that, given a $(p,q)$-form $\phi$, we have that
\begin{equation}\label{le}
\Delta'_{p''} \phi= \sum_{a=1}^l <e^a_{p-1,q},\del^*\phi> \del e^a_{p-1,q} +\sum_{b=1}^m < e_{p+1,q}^b,\del \phi>\del^* e_{p+1,q}^b,
\end{equation}
where $l=\Dim_{\C}\Ha^{p-1,q}$, $m=\Dim_{\C}\Ha^{p+1,q}$ and $\{e^\circ_{\bullet,\bullet}\}$ is an orthonormal basis of $\Ha^{\bullet,\bullet}$ as a $\C$-vector space and $<\bullet,\circ>$ is the standard $L^2$ product. Since, for every $(p,q)$-form $\psi$, $<e_{p,q}^\circ,\psi>\in \C$ the following observation is a direct consequence of (\ref{le}):
\begin{oss}\label{rmk}\label{rmk}
Locally, for every multiindex $L$, the derivative $D^L_j (\Delta'_{p''}\phi)_{\mid{U_j}}$ involves only the first derivatives of $\phi$.
\end{oss}

In the following we prove that the \textit{a priori estimate}, proved for example in \cite{kodaira2006complex}, still holds for the elliptic pseudo-differential operator with $\tilde\Delta$.
\begin{teo}\label{stima}
For every $k\in\Z$ there exists a constant $C_k$, depending on $k$, the constant of ellipticity of $\tilde{\Delta}$ and on the coefficients of $\tilde{\Delta}$, such that for every $\phi\in \Lambda^{p,q}(M)$
\begin{equation}\label{fe}
\norm{\phi}_{k+2}^2\leq C \left( \norm{\tilde{\Delta}\phi}^2_{k}+\norm{\phi}^2_k \right) 
\end{equation}
\end{teo}
\begin{proof}
We prove the theorem for $\Cinf$ function, then by using (\ref{SNF}) the general case of an arbitrary $(p,q)$-form follows straightforward.
Let $\phi$ be a $\Cinf$ function and let $\{\eta_j\}$ be a partition of unity subordinated to the finite open covering $\{U_j\}$. Denoting with $\phi_j:=\eta_j \phi$, by the definition (\ref{SN}) we have
\begin{equation}
\begin{array}{l}
\norm{\tilde{\Delta}\phi}^2_{k}=\norm{\Delta''\phi+\Delta'_{p''}\phi}^2_{k}=\\
\norm{\Delta''\phi}^2_{k}+\norm{\Delta'_{p''}\phi}^2_{k}+\sum\int_{U_j} \left[\left(D^l_j \Delta'' \phi_j\right)\overline{\left(D^l_j \Delta'_{p''}\phi_j\right)}\right]dX_j+\\
+\sum\int_{U_j} \left[\overline{\left(D^l_j \Delta'' \phi\right)}\left(D^l_j \Delta'_{p''}\phi\right)\right]dX_j.
\end{array}
\end{equation}
Since the sum of the last two terms is a real number, it can be either negative or non-negative. If it is non-negative, we can simply estimate
\begin{equation}
\norm{\tilde{\Delta}\phi}^2_{k}\geq \norm{\Delta''\phi}^2_{k}.
\end{equation}
By \cite{friedrichs1953differentiability}, we have that there exists a constant $C_k$ depending only by $k$, the constant of ellipticity of $\tilde{\Delta}$ and on the coefficients of $\tilde{\Delta}$, such that
\begin{equation}\label{qq}
\norm{\Delta''\phi}^2_{k}\geq C_k^{-1}\norm{\phi}^2_{k+2}-\norm{\phi}^2_k
\end{equation}
and hence we get the conclusion.

If the sum is negative, we proceed in the following way
$$
\begin{array}{l}
\sum\int_{U_j} \left[\left(D^l_j \Delta'' \phi_j\right)\overline{\left(D^l_j \Delta'_{p''}\phi_j\right)}\right]dX_j+\sum\int_{U_j} \left[\overline{\left(D^l_j \Delta'' \phi_j\right)}\left(D^l_j \Delta'_{p''}\phi_j\right)\right]dX_j=\\
 \\
-\left|\sum\int_{U_j} \left[\left(D^l_j \Delta'' \phi_j\right)\overline{\left(D^l_j \Delta'_{p''}\phi_j\right)}\right]dX_j+\sum\int_{U_j} \left[\overline{\left(D^l_j \Delta'' \phi_j\right)}\left(D^l_j \Delta'_{p''}\phi_j\right)\right]dX_j\right| \\
\\
\mbox{by Stokes Theorem the second row is equal to}\\
\\
-\left| \sum(-1)^l\int_{U_j} \left[\left(\phi_j\right)\left(\Delta'' D^l_j \overline{D^l_j \Delta'_{p''}\phi_j}\right)\right]dX_j+(-1)^l\int_{U_j} \left[\overline{\phi_j}\left(\overline{\Delta'' D^l_j}D^l_j \Delta'_{p''}\phi_j\right)\right]dX_j\right|\geq \\
\\
-\sum\left|\int_{U_j} \left[\phi_j\left(\Delta'' D^l_j \overline{D^l_j \Delta'_{p''}\phi_j}\right)\right]dX_j+\int_{U_j} \left[\overline{\phi_j}\left(\overline{\Delta'' D^l_j}D^l_j \Delta'_{p''}\phi_j\right)\right]dX_j\right| \geq\\
\\
-\sum\int_{U_j} \left|\left[\phi_j\left(\Delta'' D^l_j \overline{D^l_j \Delta'_{p''}\phi_j}\right)\right]+\left[\overline{\phi_j}\left(\overline{\Delta'' D^l_j}D^l_j \Delta'_{p''}\phi_j\right)\right]\right|dX_j\geq\\
\\
-2\sum\int_{U_j}\left| \phi_j\left(\Delta'' D^l_j \overline{D^l_j \Delta'_{p''}\phi_j}\right)\right|dX_j\geq\\
\\
-2\sum\left(\int_{U_j}|\phi_j|^2dX_j\right)^{\frac{1}{2}}\left(\int_{U_j}\left|(\Delta'' D^l_j \overline{D^l_j \Delta'_{p''}\phi_j}\right|^2dX_j\right)^{\frac{1}{2}}\geq\\
\\
-2\sum\left(\int_{U_j}|\phi_j|^2dX_j\right)^{\frac{1}{2}}\sum\left(\int_{U_j}\left|(\Delta'' D^l_j \overline{D^l_j \Delta'_{p''}\phi_j}\right|^2dX_j\right)^{\frac{1}{2}}.\\
\end{array}
$$
By Remark \ref{rmk}, we have that, for every $k\in \N$, the $k$-th Sobolev norm of $\Delta'_{p''}\phi$ can be estimated with the $0$-th norm of the first derivative of $\phi$. In fact, since $\phi$ is a $\Cinf$ function we have that
\begin{equation}
D^l(\Delta'_{p''} \phi)=\sum <e_{1,0}^a,\del\phi>D^l \del^* e_{1,0}^a,
\end{equation}
hence
\begin{equation}
\norm{\Delta'_{p''} \phi}_k^2=\sum \int_{U_j} |\sum_a <e_{1,0}^a,\del\phi>D^l \del^* e_{1,0}^a|^2 dX_j.
\end{equation}
Since $\{e_{1,0}^a\}$ is an orthonormal basis, then, for every $a=1,\dots, \dim H^{1,0}_{\delbar}(M)$, there exists a constant $C'_a>0$ such that $<e_{1,0}^a,\del\phi>\leq C'_a \norm{\del\phi}_0$. Let $C'$ be the supremum  of such constant and let
\begin{equation}
C'_k:=\sum \int_{U_j} |D^l \del^* e_{1,0}^a|^2 dX_j,
\end{equation}
then we have
\begin{equation}
\norm{\Delta'_{p''}\phi}_k\leq C'_k\norm{\phi'}_0.
\end{equation}
Thus we have
$$
-2\sum\left(\int_{U_j}|\phi|^2dX_j\right)^{\frac{1}{2}}\sum\left(\int_{U_j}\left|(\Delta'' D^l_j \overline{D^l_j \Delta'_{p''}\phi}\right|^2dX_j\right)^{\frac{1}{2}}\geq -2C'_k \norm{\phi}_0\norm{\phi'}_0\geq -C'_k\norm{\phi}_1^2.
$$
Summing up we have
\begin{equation}
\norm{\tilde{\Delta}\phi}^2_{k}\geq \norm{\Delta''\phi}^2_{k}+\norm{\Delta'_{p''}\phi}^2_{k}-C'_k\norm{\phi}_1^2.
\end{equation}
Using (\ref{qq}) we obtain
\begin{equation}
\norm{\tilde{\Delta}\phi}^2_{k}\geq  C_k^{-1}\norm{\phi}^2_{k+2}-\norm{\phi}^2_k-C'_k\norm{\phi}_1^2.
\end{equation}
If $k\geq 1$ then we have
\begin{equation}
\norm{\tilde{\Delta}\phi}^2_{k}\geq C_k^{-1}\norm{\phi}^2_{k+2}-(1+C'_k)\norm{\phi}^2_k
\end{equation}
which is equivalent to the thesis with $C=C_k(1+C'_k)$.

If $k=0$, we must be more precise in our estimates. In particular, by the definition of $\Delta'_{p''}$ and denoting with $\{e_{1,0}^m\}$ a basis of $\mathcal{H}^{1,0}$, we have that
\begin{equation}\label{16}
\sum \int_{U_j} \left|\Delta''\Delta'_{p''}\phi \right|^2 dX_j = \sum \int_{U_j} \left|\Delta''\left(\sum_{m} <\del\phi,e_{1,0}^m>\del^* e_{1,0}^m\right) \right|^2 dX_j.
\end{equation}
Since $<\del\phi,e_{1,0}^m>$ is a complex number, more precisely
\begin{equation}
<\del\phi,e_{1,0}^m>=<\phi,\del^* e_{1,0}^m> \leq C_m \norm{\phi}_0,
\end{equation}
we can estimate (\ref{16}) with the $0$-th norm of $\phi$ obtaining
\begin{equation}
\norm{\tilde{\Delta}\phi}^2_{0}\geq C_0^{-1}\norm{\phi}^2_{2}-(1+C''_0)\norm{\phi}^2_0,
\end{equation}
where $C''_0$ depends only by $\{e_{1,0}^m\}$.

\end{proof}

Theorem \ref{stima} provides the basis for the estimates in the following sections. Since we need to construct the Green operator of $\tilde\Delta$, we prove the following 
\begin{prop}
If $\zeta$ is different from any eigenvalue of $\tilde\Delta$, then $\tilde\Delta-\zeta Id$ is bijective.
\end{prop}
\begin{proof}{ \cite[p. 338]{kodaira2006complex} }
By construction it is obvious that $\tilde\Delta-\zeta Id$ is injective. We want to prove that $\tilde\Delta-\zeta Id$ is surjective. Let $\phi\in\Cinf(M,\C)$, then we can write $\phi=\sum_{j=1}^{\infty} b_j e_j$, where $\{e_j\}$ is an orthonormal basis of $\Cinf(M,\C)$ made of eigenfunctions of $\tilde{\Delta}$. By \cite[Lemma 7.4]{kodaira2006complex}, we have that, for every $l\in \N$,
$$
\sum_{j=1}^{\infty} |\lambda_j|^{2l} |b_j|^2<\infty,
$$
where $\lambda_j$ is the eigenvalue associated to $e_j$.
Let $a_j:= b_j/(\lambda_j-\zeta)$, then
$$
\sum_{j=1}^{\infty} |\lambda|^{2l} |a_j|^2=\sum_{j=1}^{\infty} \frac{|\lambda_j|^{2l}|b_j|^2}{|\lambda_j-\zeta|^2}<\infty
$$
for every $l\in\N$. Thus we have that there exists $\psi\in\Cinf(M,\C)$ such that $\psi=\sum_{j=1}^{\infty} a_j e_j$. By direct computation we have
$$
(\tilde\Delta-\zeta Id)\psi=\sum_{j=1}^{\infty} (\lambda_j-\zeta)a_j\ e_j=\sum_{j=1}^{\infty} b_j e_j=\phi
$$
and hence $\tilde\Delta-\zeta Id$ is surjective.
\end{proof}
\end{section}

\begin{section}{Small Deformations}\label{defo}
In this section we prove Theorem \ref{MT}; to do so we use the theory of deformations of complex structures developed by Kodaira and Spencer in \cite{kodaira1960deformations}. Our approach is similar to the one in \cite{kodaira2006complex} for elliptic self-adjoint differential operators and, when the proof of a statement does not change when we replace a generic differential operator with $\tilde\Delta$, we omit that proof.

Let $\{J_t\}$, with $t\in B\subset\C^m$, be a $\Cinf$ family of complex structures over $M$ and let $\{\omega_t\}$ be a $\Cinf$ family of Hermitian metrics on the fibers. Suppose that $0\in B$. For every $t\in B$, we denote with $\Delta''_t$, $\Delta'_{p''t}$ and $\tilde\Delta_t$ the operators described in the previous section with respect to the complex structure $J_t$. Moreover, in this section, we assume the following
\begin{hyp}
$\Dim \Ker \Delta''_t$ is a constant independent of $t$.
\end{hyp}
\begin{prop}
If  $\Dim \Ker \Delta''_t$ is independent of $t$, then $\{\tilde\Delta_t\}$ is a $\Cinf$ family of pseudo differential operators.
\end{prop}

\begin{proof}
From \cite{kodaira2006complex}, we have that all the derivative operators depend $\Cinf$ with respect to $t$, hence we only need to prove that if $\{\phi_t\}$ and $\{\psi_t\}$ are $\Cinf$ family of $(p,q)$ forms over $M$ and if $\{g_t\}$ is a $\Cinf$ family of Hermitian metrics over $M$, then the scalar product $<\psi_t,\phi_t>_t$ varies in a $\Cinf$ way respect to $t$.
By definition
\begin{equation}\label{ps}
<\psi_t,\phi_t>_t=\int_{M} \psi_t\wedge *_t \overline{\phi}_t.
\end{equation}

Let $\{U_j\}$ be a finite covering of $M$ made by open coordinate neighborhood and let $\{\eta_j\}$ be a partition of unity subordinate to $\{U_j\}$. Then, for every $t\in B$ we have
\begin{equation}
\int_{M} \psi_t\wedge *_t\overline{\phi}_t=\sum_{j} \int_{U_j}\eta_j \psi_t\wedge *_t\overline{\phi}_t.
\end{equation}

Now, locally we have $\psi_t=\sum_{A_p,B_q} \psi_t^{A_p\overline{B}_q} dz^{A_p\overline{B}_q}$ and 
\begin{equation}
*_t\overline\phi_t(z):= (i)^n(-1)^k \sum_{A_p,B_q} sgn\left(\begin{array}{cc} A_p & A_{n-p}\\ B_q & B_{n-q} \end{array}\right) g_t(z)\overline\phi_t^{A_p\overline{B}_q}(z) dz^{B_{n-p}\overline{A_{n-q}}}.
\end{equation}
Then we have
\begin{equation}
\int_{U_j}\eta_j \psi_t\wedge *_t\overline{\phi}_t=\sum_{A_p,B_q} \sigma_{A_p B_p}\int_{U_j} \eta_j \psi_t^{A_p\overline{B}_q}\overline{\phi}_t^{A_p\overline{B}_q}g_t dz^1\wedge \dots \wedge dz^{n}\wedge d\overline{z}^1 \wedge \dots \wedge d\overline{z}^n,
\end{equation}
where $\sigma_{A_p B_p}$ is the sign of the permutation.

In order to prove that the scalar product is $\Cinf$, it suffices to show that it is $\mathcal{C}^k$ for every $k \in \N$. We begin proving that it is $\mathcal{C}^0$: by the continuity of the integral operator, the coefficients of $\psi_t$, $\phi_t$ and $g_t$ and since $\eta_j \psi_t\wedge *_t\overline{\phi}$ is continuous and compactly supported in $U_j$ , we have  the continuity of the scalar product.

Now we prove by induction over $r\in \N$ that (\ref{ps}) is $\mathcal{C}^r$. Let $r=1$ and consider the following
\begin{equation}
\frac{<\psi_t,\phi_t>_t-<\psi_s,\phi_s>_s}{t-s}=\frac{1}{t-s}\int_M \psi_t\wedge *_t \overline{\phi}_t-\psi_s\wedge*_s\overline{\phi}_s.
\end{equation}
Locally we can rewrite the integral above as
\begin{equation}
\frac{1}{t-s}\sum_{A_p,B_q} \int_{U_j} \left( \eta_j \psi_t^{A_p\overline{B}_q}\overline{\phi}_t^{A_p\overline{B}_q}g_t - \eta_j \psi_s^{A_p\overline{B}_q}\overline{\phi}_s^{A_p\overline{B}_q}g_s \right)dz^1\wedge \dots \wedge dz^{n}\wedge d\overline{z}^1 \wedge \dots \wedge d\overline{z}^n.
\end{equation}
Now, for every multi-indexes $A_p$ and $B_q$, we consider the following construction
\begin{equation}
\begin{array}{ll}
\frac{1}{t-s}\eta_j\left( \psi_t^{A_p\overline{B}_q}\overline{\phi}_t^{A_p\overline{B}_q}g_t - \psi_s^{A_p\overline{B}_q}\overline{\phi}_s^{A_p\overline{B}_q}g_s\right)= & \\
 = \frac{1}{t-s}\eta_j\left( \psi_t^{A_p\overline{B}_q}\overline{\phi}_t^{A_p\overline{B}_q}(g_t-g_s) + \psi_t^{A_p\overline{B}_q}(\overline{\phi}_t^{A_p\overline{B}_q}-\overline{\phi}_s^{A_p\overline{B}_q})g_s+ (\psi_t^{A_p\overline{B}_q}-\psi_s^{A_p\overline{B}_q})\overline{\phi}_s^{A_p\overline{B}_q}g_s\right) & \\
  = \eta_j\left( \psi_t^{A_p\overline{B}_q}\overline{\phi}_t^{A_p\overline{B}_q}\frac{g_t-g_s}{t-s} + \psi_t^{A_p\overline{B}_q}\frac{\overline{\phi}_t^{A_p\overline{B}_q}-\overline{\phi}_s^{A_p\overline{B}_q}}{t-s}g_s+ \frac{\psi_t^{A_p\overline{B}_q}-\psi_s^{A_p\overline{B}_q}}{t-s}\overline{\phi}_s^{A_p\overline{B}_q}g_s\right). &
\end{array}
\end{equation}
When $t$ tends to $s$, we obtain that 
\begin{equation}\label{der}
\begin{array}{rl}
\frac{d <\psi_t,\phi_t>_t}{dt}\mid_{t=s}= & <\psi'_s,\phi_s>_s+<\psi_s,\phi'_s>_s \\
+ & \sum_{j} \sum_{A_p,B_q} \int_{U_j} \eta_j \psi_s^{A_p\overline{B}_q}\overline{\phi}_s^{A_p\overline{B}_q}g'_s dz^1\wedge \dots \wedge dz^{n}\wedge d\overline{z}^1 \wedge \dots \wedge d\overline{z}^n,
\end{array}
\end{equation}
where $\phi'_s$ and $\psi'_s$ denote the derivative along $t$ of the $\Cinf$ forms $\phi(z,t)$ and $\psi(z,t)$ respectively.  
Using the same arguments of the $\mathcal{C}^0$ case, we have the derivative is continuous.

Suppose, by induction, that (\ref{ps}) is $\mathcal{C}^r$. By reiteration of (\ref{der}) we have that the $r$-th derivative of (\ref{ps}) is made by two type of components
\begin{itemize}
\item[i)] $<\psi_t^{(i)},\phi_t^{(j)}>_t$;
\item[ii)] $\int_M g^{(k)}_t\psi_t^{(i)}\wedge\phi_t^{(j)}$, for some $k\in \N$.
\end{itemize}
In either case, using the same argument as above, we have the existence and the continuity of the derivative of the $r$-th derivative of $<\phi(z,t),\psi(z,t)>_t$.


\end{proof}

\setcounter{teo}{2}
\begin{teo}
If  $\Dim \Ker \Delta''_t$ is independent of $t$, then $Dim \Ker \tilde\Delta_t$ is an upper semi-continuous function in $t$.
\end{teo}
\setcounter{teo}{6}
In order to prove Theorem \ref{ker} we need some preliminary results.
 
\begin{lemma}[Fridedichs' Inequality]\label{fri}
For every $k\in \N$, there exists a constant $C_k$ independent of $t$ such that, for every $\Cinf$ family $\{\phi_t\}$ with $\phi_t\in\Lambda^{p,q}(M,J_t)$, the inequality
\begin{equation}\label{eq}
\norm{\phi}_{k+2}^2\leq C_k(\norm{\tilde\Delta_t\phi}_k^2+\norm{\phi}_0^2)
\end{equation}
holds.
\end{lemma}

\begin{proof}
By Theorem \ref{stima}, for every $t\in B$ and for every $k\in \N$ there exists a constant $C_{k,t}$ such that
\begin{equation}
\norm{\phi_t}_{k+2}^2\leq C_{k,t}(\norm{\tilde\Delta_t\phi_t}_k^2+\norm{\phi_t}_k^2)
\end{equation}
holds for every $\phi\in \Lambda^{p,q}(M,J_t)$.
Since $\{g_t\}$ is a continuous family the Sobolev norm varies continuously with respect to $t$, then, up to shrinking $B$, we have that $C_k:=sup\{C_{k,t}\}<\infty$. Then we obtain
\begin{equation}\label{eq1}
\norm{\phi_t}_{k+2}^2\leq C_{k}(\norm{\tilde\Delta_t\phi_t}_k^2+\norm{\phi_t}_k^2).
\end{equation}

We will prove (\ref{eq}) by induction over $k$. For $k=0$, the equations (\ref{eq}) and (\ref{eq1}) are the same. Then the thesis holds for $k=0$. Now suppose that thesis holds for $k-1$. By equation (\ref{eq}) at the step $k-1$ we have
$$
\norm{\phi}_k^2\leq \norm{\phi}_{k+1}^2\leq C_{k-1}\left( \norm{\tilde\Delta\phi}^2_{k-1} +\norm{\phi}_0^2\right).
$$
Since $\norm{\tilde\Delta\phi}^2_{k-1}\leq \norm{\tilde\Delta\phi}^2_{k}$, we get
\begin{equation}
\norm{\phi}_{k+2}^2\leq C_k(\norm{\tilde\Delta_t\phi}_k^2+\norm{\phi}_k^2)\leq C_k(1+C_{k-1})\left(\norm{\tilde\Delta_t\phi}_k^2+\norm{\phi}_0^2\right).
\end{equation}

\end{proof}

\begin{teo}
If $E_t$ is bijective for every $t$ and if there exists a constant $c$ such that $\norm{\phi}_0\leq c\norm{E_t \phi_t}_0$ for every $\phi_t\in\Lambda^{p,q}(M,J_t)$. Then the family $\{G_t\}$ of Green operator associated to $\{E_t\}$ is $\Cinf$ differentiable in $t$.
\end{teo}
\begin{proof}
This theorem can be proved using the same arguments of \cite[Theorem 7.5]{kodaira2006complex} by making the following observation: by hypothesis, we have that, for every $t$
\begin{equation}
\norm{\phi_t}_0\leq c\norm{E_t \phi_t}_0\leq c\norm{E_t \phi_t}_k.
\end{equation}
Hence, there exists a constant $c'$ independent of $t$ and $\phi_t$ such that
\begin{equation}
\norm{\phi_t}_{k+2}\leq c'\norm{\tilde\Delta_t \phi_t}_k.
\end{equation}
By Sobolev's inequality, we have for large enough $k$:
\begin{equation}
|D^l_j\phi_j^{\lambda}(x)|\leq c\norm{\phi}_{k+2}.
\end{equation}
Hence
\begin{equation}
 |D^l_j\phi_j^{\lambda}(x)|\leq C\norm{\tilde{\Delta}_t\phi}_{k}.
\end{equation}

Now we are in the same conditions of \cite[Theorem 7.5]{kodaira2006complex} and we can proceed in same way.


\end{proof}
In order to proceed we need the following result that guarantees the existence of an orthonormal basis of eigenvectors.

\begin{lemma}[{\cite[Lemma 1.6.3]{gilkey1984invariance}}]
Let $P:\Cinf(M,\C)\rightarrow \Cinf(M,V)$ be an elliptic self-adjoint pseudo-differential operator of order d> 0. Then
\begin{itemize}
\item We can find a complete orthonormal basis $\{\psi_n\}_{n=1}^{\infty}$ for $L^2(M)$ of eigenvectors of $P$. $P\psi_n=\lambda_n\psi_n$.
\item The eigenvectors $\psi_n$ are smooth and $\lim_{n\rightarrow\infty} |\lambda_n| = \infty$.
\item If we order the eigenvalues $|\lambda_1|\leq|\lambda_2|\leq\dots$ then there exists a constant $C > 0$ and an exponent $\delta> 0$ such that $|\lambda_n| \geq Cn\delta$ if $n>n_0$ is large.
\end{itemize}
\end{lemma}

\begin{lemma}
Let $\zeta_0\in\C$ be different from any eigenvalue of $\tilde\Delta_0$. Then there exist $\delta>0$ and $c>0$ such that, for $|t|<\delta$ and $|\zeta-\zeta_0|<0$, the following inequality
\begin{equation}
\norm{\phi}_0\leq c\norm{(\tilde\Delta_t-\zeta Id)\phi}_0
\end{equation}
holds.
\end{lemma}
\begin{proof}
Suppose that, for any small $\delta>0$ there is no such constant. Then, for q=1,2,..., there exist $t_q\in B$, $\zeta_q\in \C$ and $\phi_q \in \Lambda^{p,q}(M_t)$ such that
$$
\begin{array}{cccc}
|t_q-t_0|<\frac{1}{q}, & |\zeta_q-\zeta_0|<\frac{1}{q}, & \norm{\tilde{\Delta}_{t_q}(\zeta_q)\phi_q}_0<\frac{1}{q}, & \norm{\phi_q}_0=1.
\end{array}
$$
By Lemma \ref{fri}, we have
\begin{equation}\label{ooo}
\norm{\phi_q}^2_2\leq c_0(\norm{\tilde{\Delta}_{t_q}(\zeta_q)\phi_q}_0^2+\norm{\phi_q}_0^2)\leq 2c_0.
\end{equation}
Since the coefficients of $\tilde{\Delta}$ are $\mathcal{C}^\infty$ function on $(x,t)$ and $\tilde{\Delta}(\zeta_q)\phi_q$ is uniformly bounded, we have
\begin{equation}
\norm{\tilde{\Delta}_{t_q}(\zeta_q)\phi_q-\tilde{\Delta}_{0}(\zeta_0)\phi_q}_0\rightarrow 0.
\end{equation}
Since $\norm{\tilde{\Delta}_{t_q}(\zeta_q)\phi_q}_0<\frac{1}{q}$ then
$$
\norm{\tilde{\Delta}_0(\zeta_0)\phi_q}_0\rightarrow 0.
$$
But by (\ref{ooo}) for a suitable $\mu_0$ we have
$$
\norm{\tilde{\Delta}_0(\zeta_0)\phi_q}_0\geq \mu_0 \norm{\phi_q}_0.
$$
Hence $\norm{\phi_q}_0\rightarrow 0$, which contradicts the hypothesis.
\end{proof}

Let $C$ be a closed Jordan curve on the complex which does not pass through any eigenvalue of $\tilde{\Delta}_0$. As in \cite{kodaira2006complex}, we denote with $((C))$ the interior of $C$. Since, for sufficient small $t$, $C$ does not contain any eigenvalue of $\tilde{\Delta}$, we can define the operator
\begin{equation}
F_t(C) \phi:= \sum_ {\lambda_k(t) \in ((C))} <\phi,e^k_t> e^k_t,
\end{equation} 
and put $\mathbb{F}_t(C):=F_t(C)\Lambda^{p,q}_t(M)$. Obviously, $\mathbb{F}_t(C)$ is a finite dimensional linear subspace of $\Lambda^{p,q}_t(M)$.
\begin{lemma}
The operator $F_t(C)$ can be written as
$$
F_t(C)\phi=-\frac{1}{2\pi i}\int_C G_t(\zeta)\phi d\zeta,
$$
where $G_t(\zeta)$ is the Green operator associated to $\tilde\Delta_t(\zeta)$.
\end{lemma}

\begin{lemma}
$F_t(C)$ is $\mathcal{C}^{\infty}$ in $t$ for $|t|<\delta$.
\end{lemma}
\begin{proof}
Since $G(\zeta)$ is $\mathcal{C}^\infty$ in $(t,\zeta)$, it follows that, if $\phi_t$ is $\mathcal{C}^\infty$ in $t$, then, by the definition above, also $F_t(C)$ is $\mathcal{C}^\infty$ in $t$.
\end{proof}

\begin{lemma}\label{lemma5}
$\Dim \mathbb{F}_t(C)$ is independent of $t$ for $|t|<\delta$.
\end{lemma}
\begin{proof}
Let $d=\Dim \mathbb{F}_0(C)$ and let $\{e_1,...,e_d\}$ be a basis for $\mathbb{F}_0(C)$. Since $F_t(C)e_r$ are $\mathcal{C}^\infty$ differentiable in $t$, and $F_0(C)e_r=e_r$ are linearly independent. Therefore $F_t(C)e_r$ are linearly independent for sufficient small $t$. Hence $\Dim \mathbb{F}_t\geq d$.

Suppose that, for any $\delta>0$, there exists $t$ such that $|t|<\delta$ and $\Dim F_t(C)>d$. Then it is possible to find a sequence $t_q$ such that $|t_q|<\frac{1}{q}$ and $\Dim \mathbb{F}_{t_q}>d$. Let $\{\lambda^q_r\}$ be $d+1$ eigenvalues of $\tilde{\Delta}_{t_q}$ and let $\{e^q_r\}$ be the relative eigenfunctions. For a sufficient large $k$  we have
$$
|p^q_r D_j^j e^q_r(x)|^2\leq |D_j^j e^q_r(x)|^2\leq C_k(1+\sum_{\alpha=1}^k |\lambda^q_r|^2).
$$
Then 
$$
|\tilde{\Delta}_{t_q} e^q_r(x)|^2\leq 4C_k^2(1+\sum_{\alpha=1}^k |\lambda^q_r|^2).
$$
Since $\lambda^q_r\in ((C))$, the sequence $\tilde{\Delta}_{t_q} e^q_r(x)$ is uniformly bounded in $B$. Then, up to subsequences, $\{\tilde{\Delta}_{t_q} e^q_r(x)\}$ converges uniformly  and we have
$$
\lim_{q\rightarrow +\infty} \tilde\Delta_{t_q} e^q_r =\tilde\Delta_0 e^0_r=\tilde\Delta e_r.
$$
Since $\tilde{\Delta}_{t_q} e^q_r=\lambda^q_r e^q_r$, the sequence $\{\lambda_r^q\}$ converges to $\lambda_r$ and we have
$$
\begin{array}{cc}
\tilde{\Delta}_0 e_r=\lambda_r e_r & \norm{e_r}_0=1.
\end{array}
$$
This means that $\lambda_r$ is an eigenvalue and $e_r$ an eigenfunction of $\tilde\Delta_0$. Since there are no eigenvalue on $C$, $\lambda_r\in ((C))$. Moreover
$$
(e_r,e_s)=\lim_{q\rightarrow +\infty} (e_r^q,e_s^q)=\delta_{rs},
$$
this means that there are d+1 linearly independent eigenfunctions. This contradicts the hypothesis.
\end{proof}

\begin{proof}[Proof of Theorem \ref{ker}]
To conclude the proof we only need to consider a Jordan curve passing around the eigenvalue $0$ and not containing any other eigenvalue of $\tilde\Delta_0$. The thesis follows directly from Lemma \ref{lemma5}.

\end{proof}

\setcounter{teo}{0}
\begin{teo}
Let $(M,J_t)$ be a family of complex manifolds and suppose that the dimension of $\Ker\Delta''_t\cap\Lambda^{p,q}(M,J_t)$ is independent of $t$ for every $(p,q)\in\Z^2$. Then the degeneration at the second step of the Fr\"olicher spectral sequence is stable under small deformations of the complex structure.
\end{teo}
\begin{proof}
We denote with $b_k$ the dimension of $H^k_{dR}(M;\C)$, with $\tilde{h}^{p,q}_t$ the complex dimension of $\Ker\tilde\Delta_t\cap \Lambda^{p,q}(M,J_t)$ and with $e^{p,q}_{2}(t)$ the dimension of $E^{p,q}_{2}(M,J_t)$. We recall the degeneration at the second step of $\{E^{p,q}_r(M,J_t)\}$ is equivalent to
\begin{equation}
b_k= \sum_{p+q=k} e^{p,q}_{2}(t);
\end{equation}
by Theorem \cite[Theorem 3.4]{popovici2016degeneration} we have
\begin{equation}
\sum_{p+q=k} \tilde{h}^{p,q}_t= \sum_{p+q=k} e^{p,q}_{2}(t);
\end{equation}
finally, by Theorem \ref{ker}, we know that $\tilde{h}^{p,q}_t$ is an upper-semi continuous function of $t$.
 
Suppose that the Fr\"olicher spectral sequence of $(M,J_0)$ degenerates at the second step, then, summing up all the previous considerations, we have
\begin{equation}
b_k=\sum_{p+q=k} e^{p,q}_{2,0}=\sum_{p+q=k} \tilde{h}^{p,q}_0\geq \sum_{p+q=k} \tilde{h}^{p,q}_t= \sum_{p+q=k} e^{p,q}_{2,t}\geq b_k.
\end{equation}
Thus
\begin{equation}
b_k= \sum_{p+q=k} \tilde{h}^{p,q}_t= \sum_{p+q=k} e^{p,q}_{2,t},
\end{equation}
that means that, for $t$ small enough, the Fr\"olicher spectral sequence of $M_t$ degenerates at the second step.
\end{proof}
\end{section}

\begin{section}{Example}
In this section we provide an example of $\Cinf$ curve of compact complex manifolds such that the Fr\"olicher spectral sequence degenerate at the second step for one of them and at higher steps for for the others. In this example the dimension of $\Ker \Delta''_t$ is not constant with respect to $t$ showing the importance of Assumption \ref{ass}. Let $X$ be the Nakamura manifold, namely a compact complex three-dimensional holomorphically parallelizable solvmanifold constructed in the following way: let $G$ be the Lie group defined as $G:=\C\ltimes_{\phi} \C^2$, where
\begin{equation}
\phi(z):=
\left(
\begin{array}{cc}
e^z & 0\\
0 & e^{-z}
\end{array}
\right).
\end{equation}
Let $\Gamma:=\Gamma'\ltimes_{\phi}\Gamma''$ be a lattice in $G$, where $\Gamma''$ is a lattice in $\C^2$  and $\Gamma':= \Z(a+ib)+\Z(c+id)$ is such that it is a lattice in $\C$ and $\phi(a+ib)$ and $\phi(c+id)$ are conjugate elements in $SL(4,\Z)$. Then $X:=\Gamma\setminus G$.

We consider the following deformation of $X$: let $t\in \C$ and consider the $(0,1)$-form on $X$ with value in $T^1,0X$ defined by
$$
\phi_t=te^{z_1}d\overline{z}_1\otimes\frac{\partial}{\partial z_2}.
$$
For $|t|<\epsilon$, let $X_t$ be the small deformation of $X$ associated to $\phi_t$.   
We prove the following
\setcounter{teo}{8}
\begin{teo}\label{ex}
The Fr\"olicher spectral sequence of $X_t$ degenerates at the second step for $t=0$, while it degenerates at higher steps for $t\neq 0$.
\end{teo}
First we recall that the Betti's number of $X$ are the following
\begin{equation}
\begin{array}{cccc}
b_0=b_6=1, & b_1=b_5=2, & b_2=b_4=5, & b_3=8. 
\end{array}
\end{equation}
Moreover it is a well-known fact that they do not change under deformations of the complex structure. Now we explicitly compute $E^{\bullet,\bullet}_1(X_t)$ and $E^{\bullet,\bullet}_2(X_t)$ and we show that 
\begin{equation}
b_k=\sum_{p+q=k} \Dim E^{p,q}_2(X_0)
\end{equation}
for every $k=1,...,6$, while the equality is false for $t\neq 0$.

We proceed in the following way: we begin with the computation of the Dolbeault cohomology of $X_t$ since, as we recalled is Section \ref{pre}, the first step of the Fr\"olicher spectral sequence is isomorphic to the Dolbeault cohomology, namely $E^{p,q}_1(X_t)\simeq H^{p,q}_{\delbar}(X_t)$ for every $(p,q)\in \Z^2$. By applying \cite[Theorem 1.3]{angella2012bott}, in \cite{tomassini2014dolbeault} Tomassini and Torelli found the $\Delta''_t$-harmonic forms of $X_t$; for every $(p,q)\in\Z$, those forms are a basis for $H^{p,q}_{\delbar}(X_t)$ as $\C$ vector space. Then, as proved in \cite{cordero1997general}, $E^{p,q}_2(X_t)$ can be described as
\begin{equation}
E^{p,q}_2=\frac{X^{p,q}_2(X_t)}{Y^{p,q}_2(X_t)},
\end{equation}
where
\begin{equation}
X^{p,q}_2(X_t):=\left\{\alpha \in\Lambda^{p,q}(X_t)| \delbar\alpha=0 \mbox{ and }\exists \beta\in \Lambda^{p+1,q-1}(X_0) \mbox{ s.t. } \del\alpha+\delbar\beta=0\right\},
\end{equation}
\begin{equation}
Y^{p,q}_2(X_t):=\left\{\del\alpha+\delbar\beta\in\Lambda^{p,q}(X_0) | \delbar\alpha=0\right\}.
\end{equation}
Namely, if $\alpha\in X^{p,q}_2(X_t)$, $[\alpha]_2\in E^{p,q}_2(X_t)$ can be written as
\begin{equation}
[\alpha]_2=\{\alpha+\del\beta+\delbar\gamma | \delbar\beta=0\}.
\end{equation}
On the other hand,  $E^{p,q}_2(X_0)$ is the cohomology group of the complex
\begin{equation}
E^{p-1,q}_1(X_0)\xrightarrow{d_1} E^{p,q}_1(X_0)\xrightarrow{d_1} E^{p+1,q}_1(X_t)
\end{equation}
where $d_1$ is the operator $[\alpha]\mapsto [\del\alpha]$. Since it is well defined and every $\delbar$-closed form $\alpha$ belongs to a class in the Dolbeault cohomology, we need only to work with $\Delta''_t$-harmonic form. In fact if an harmonic form $\phi$ is such that $\del\phi$ is $\delbar$-exact, then every form $\alpha=\phi+\delbar\beta$ is such that $\del\alpha$ is $\delbar$-exact. Moreover, from the decomposition
\begin{equation}
\Ker \delbar_t=\mathcal{H}^{\bullet,\bullet}(X_t)\oplus \Ima\delbar_t
\end{equation}
we have that if $\phi$ is $\del$-exact then it is the the image through $\del$ of another harmonic form.

Now we proceed with the computation.
We recall that in \cite{tomassini2014dolbeault}, the Dolbeault cohomology of $X_t$ was computed using a suitable sub-complex $(B_t,\delbar_t)\subset(\Gamma^{\bullet,\bullet}X_t,\delbar_t)$. Namely, let
$$
B_t:=\wedge^{\bullet,\bullet}(\C<\phi^{1,0}_1(t), \phi^{1,0}_2(t), \phi^{1,0}_3(t)>\oplus \C<\phi^{0,1}_1(t), \phi^{0,1}_2(t), \phi^{0,1}_3(t)>,
$$
where
$$\begin{array}{ccc}
\phi^{1,0}_1(t):=dz_1, & \phi^{1,0}_2(t):=e^{-z_1}dz_2-td\overline{z}_1, & \phi^{1,0}_3(t):=e^{z_1}dz_3,\\
\phi^{0,1}_1(t):=d\overline{z}_1, & \phi^{1,0}_2(t):=\overline{t}e^{\overline{z}_1-z_1}dz_1-e^{-z_1}d\overline{z}_2, & \phi^{1,0}_2(t):=e^{z_1}d\overline{z}_3.
\end{array}
$$
We have the following  structure equations
$$
\left\{
\begin{array}{l}
d\phi^{1,0}_1(t)=0\\
d\phi^{1,0}_2(t)=-\phi^{1,0}_1(t)\wedge\phi^{1,0}_2(t)-t\phi^{1,0}_1(t)\wedge\phi^{0,1}_1(t)\\
d\phi^{1,0}_3(t)=\phi^{1,0}_1(t)\wedge\phi^{1,0}_3(t)\\
d\phi^{0,1}_1(t)=0\\
d\phi^{0,1}_2(t)=-\phi^{1,0}_1(t)\wedge\phi^{0,1}_2(t)+\overline{t}e^{\overline{z}_1-z_1}\phi^{1,0}_1(t)\wedge\phi^{0,1}_1(t)\\
d\phi^{0,1}_3(t)=\phi^{1,0}_1(t)\wedge\phi^{0,1}_3(t)\\
\end{array}
\right.
$$
Then $(B_t,\delbar_t)$ is a finite sub-complex of $(\Gamma^{\bullet,\bullet}X_t,\delbar_t)$, which is smooth on $X\times B(0,\epsilon)$, closed with respect to the $\C$-anti-linear Hodge star operator $*_t$ associated to the Hermitian metric $g_t:= \sum_{i=1}^3 \phi_i^{1,0}(t) \odot \overline{\phi_i^{1,0}(t)}$ and such that $H^{\bullet,\bullet}_{\delbar}(B_0)\simeq H^{\bullet,\bullet}_{\delbar}(X)$.

We start with the computation of $E^{p,q}_1(X_0)$. For $t=0$ we have that every $(p,q)$-form of the type
$$
\wedge^{p,q}\left( \C< \phi^{1,0}_1(0),\phi^{1,0}_2(0),\phi^{1,0}_3(0)>\oplus \C< \phi^{0,1}_1(0),\phi^{0,1}_2(0),\phi^{0,1}_3(0)>\right)
$$
is $\Delta''$-harmonic. Thus they form a basis for $E^{p,q}_1(X_0)$. Next we compute the image through $\del$ of every harmonic form. In Table \ref{tab1} we report only the harmonic forms that are not $\del$-closed.

 Using the same Table, we also know which harmonic form is $\del$-exact and so represents the zero class in $E^{p,q}_2(X_0)$. Let $\phi$ be a harmonic $(p,q)$-form in Table \ref{tab1}, we want to know if $\del\phi$ is $\delbar$-exact. To do so we consider the scalar product of $\del\phi$ with every $(p+1,q)$ harmonic form and, by straightforward computation, we have that it is never zero. Then, again by the decomposition $\Ker \delbar_t= \mathcal{H}^{\bullet,\bullet}_{\delbar}(X_0)\oplus \Ima \delbar_t$,  they are not such that $\del\phi$ is $\delbar$-exact. In Table \ref{tab2} is listed a basis for $E^{p,q}_2(X_0)$.


%

Instead, for $t\neq 0$, the $\Delta''$-harmonic forms are listed in Table \ref{tab3}

Using the same arguments as in the case $t=0$, we obtain that $E^{p,q}_2(X_t)$ is generated by the forms listed in Table \ref{tab4}

Let $h^{p,q}(t):=\Dim H^{p,q}_{\delbar}(X_t)$ and $e^{p,q}_2(t):=\Dim E^{p,q}_2(X_t)$.

 Then we have that, for $t=0$ and $k\in\Z$, $\sum_{p+q=k} e^{p,q}_2(0)=b_k$ and this is equivalent to the degeneration at the second step of $(E^{\bullet,\bullet}_r(X_0),d_r)$. While, for $t\neq 0$ and $k=2,4$, $b_k<\sum_{p+q=k} e^{p,q}_2(t)$ (see Table \ref{tab5}).

So we have proved Theorem \ref{ex}.
\begin{oss}
The dimensions of $E^{p,q}_2(X_t)$ do not vary, in general, nor upper semi-continuously neither lower semi-continuously with respect to $t$. For example $e^{2,1}_2(t)$ is a lower semi-continuous function of $t$ and $e^{2,2}_2(t)$ is upper semi-continuous.
\end{oss}

\end{section}

\bibliographystyle{alpha}
\bibliography{biblio}

\begin{table}[!h]
\begin{tabular}{l|l}
\hline
$\phi$ & $\del\phi$\\
\hline
 & \\
$\phi^{1,0}_2$ & $-\phi^{1,0}_{12}$ \\
& \\
$\phi^{1,0}_3$ & $\phi^{1,0}_{13}$\\
 & \\
$\phi^{0,1}_{2}$ & $-\phi^{1,0}_{1}\wedge \phi^{0,1}_{2}$ \\
 & \\
$\phi^{0,1}_{3}$ & $\phi^{1,0}_{1}\wedge \phi^{0,1}_{3}$ \\
 & \\
$\phi^{1,0}_{2}\wedge\phi^{0,1}_{1}$ & $-\phi^{1,0}_{12}\wedge\phi^{0,1}_{1}$\\
 & \\
$\phi^{1,0}_{2}\wedge\phi^{0,1}_{2}$ & $-2\phi^{1,0}_{12}\wedge\phi^{0,1}_{2}$\\
 & \\
$\phi^{1,0}_{3}\wedge\phi^{0,1}_{1}$ & $\phi^{1,0}_{13}\wedge\phi^{0,1}_{1}$\\
 & \\
$\phi^{1,0}_{3}\wedge\phi^{0,1}_{3}$ & $\phi^{1,0}_{13}\wedge\phi^{0,1}_{3}$\\
 & \\
$\phi^{1,0}_{23}\wedge\phi^{0,1}_{2}$ & $-\phi^{1,0}_{123}\wedge\phi^{0,1}_{2}$\\
 & \\
$\phi^{1,0}_{23}\wedge\phi^{0,1}_{3}$ & $\phi^{1,0}_{123}\wedge\phi^{0,1}_{3}$\\
 & \\
$\phi^{0,1}_{12}$ & $-\phi^{1,0}_{1}\wedge\phi^{0,1}_{12}$\\
 & \\
$\del\phi^{0,1}_{13}$ & $\phi^{1,0}_{1}\wedge\phi^{0,1}_{13}$\\
 & \\
$\phi^{1,0}_{2}\wedge\phi^{0,1}_{12}$ & $-2\phi^{1,0}_{12}\wedge\phi^{0,1}_{12}$\\
 & \\
$\phi^{1,0}_{2}\wedge\phi^{0,1}_{23}$ & $-\phi^{1,0}_{12}\wedge\phi^{0,1}_{23}$\\
 & \\
$\phi^{1,0}_{3}\wedge\phi^{0,1}_{13}$ & $2\phi^{1,0}_{13}\wedge\phi^{0,1}_{13}$\\
 & \\
$\phi^{1,0}_{3}\wedge\phi^{0,1}_{23}$  & $\phi^{1,0}_{13}\wedge\phi^{0,1}_{23}$\\
 & \\
$\phi^{1,0}_{23}\wedge\phi^{0,1}_{12}$ & $-\phi^{1,0}_{123}\wedge\phi^{0,1}_{12}$\\
 & \\
$\phi^{1,0}_{23}\wedge\phi^{0,1}_{13}$ & $\phi^{1,0}_{123}\wedge\phi^{0,1}_{13}$\\
 & \\
$\phi^{1,0}_{2}\wedge\phi^{0,1}_{123}$ & $-\phi^{1,0}_{12}\wedge\phi^{0,1}_{123}$\\
 & \\
$\phi^{1,0}_{3}\wedge\phi^{0,1}_{123}$ & $\phi^{1,0}_{13}\wedge\phi^{0,1}_{123}$\\
\hline
\end{tabular}
\caption{Image through $\del$ of $\Delta''$-harmonic forms of $X_0$.}\label{tab1}
\end{table}

\begin{table}[!th]
\begin{tabular}{l|c}
\hline
$(p,q)$ & Basis for $E^{p,q}_2(X_0)$\\
\hline
 & \\
$(0,0)$ & $1$\\
 & \\
$(1,0)$ & $\phi^{1,0}_1$ \\
 & \\
$(0,1)$ & $\phi^{0,1}_1$\\
 & \\
$(2,0)$ & $\phi^{1,0}_{23}$\\
 & \\
$(1,1)$ & $\phi^{1,0}_{1}\wedge\phi^{0,1}_{1}$, $\phi^{1,0}_{2}\wedge\phi^{0,1}_{3}$, $\phi^{1,0}_{3}\wedge\phi^{0,1}_{2}$\\
 & \\
$(0,2)$ & $\phi^{0,1}_{23}$\\
 & \\
$(3,0)$ & $\phi^{1,0}_{123}$\\
 & \\
$(2,1)$ & $\phi^{1,0}_{12}\wedge\phi^{0,1}_{3}$, $\phi^{1,0}_{13}\wedge\phi^{0,1}_{2}$, $\phi^{1,0}_{23}\wedge\phi^{0,1}_{1}$\\
 & \\
$(1,2)$ & $\phi^{1,0}_{1}\wedge\phi^{0,1}_{23}$, $\phi^{1,0}_{2}\wedge\phi^{0,1}_{13}$, $\phi^{1,0}_{3}\wedge\phi^{0,1}_{12}$\\
 & \\
$(0,3)$ & $\phi^{0,1}_{123}$\\
 & \\
$(3,1)$ & $\phi^{1,0}_{123}\wedge\phi^{0,1}_{1}$\\
 & \\
$(2,2)$ & $\phi^{1,0}_{12}\wedge\phi^{0,1}_{13}$, $\phi^{1,0}_{13}\wedge\phi^{0,1}_{12}$, $\phi^{1,0}_{23}\wedge\phi^{0,1}_{23}$\\
 & \\
$(1,3)$ & $\phi^{1,0}_{1}\wedge\phi^{0,1}_{123}$\\
 & \\
$(3,2)$ & $\phi^{1,0}_{123}\wedge\phi^{0,1}_{23}$\\
 & \\
$(2,3)$ & $\phi^{1,0}_{23}\wedge\phi^{0,1}_{123}$\\
 & \\
$(3,3)$ & $\phi^{1,0}_{123}\wedge\phi^{0,1}_{123}$\\
\hline
\end{tabular}
\caption{ Basis for $E^{p,q}_2(X_0)$.}\label{tab2}
\end{table}

\begin{table}
\begin{tabular}{l|l}
\hline
Bi-degree & $\Delta''$-Harmonic form\\
\hline
 & \\
(0,0) & 1\\
 & \\
(1,0) & $\phi^{1,0}_{1}(t)$, $\phi^{1,0}_{3}(t)$\\
 & \\
(2,0) & $\phi^{1,0}_{12}(t)$, $\phi^{1,0}_{13}(t)$\\
 & \\
(3,0) & $\phi^{1,0}_{123}(t)$ \\ 
 & \\
(0,1) & $\phi^{0,1}_{1}(t)$, $\phi^{0,1}_{2}(t)$, $\phi^{0,1}_{3}(t)$\\
 & \\
(1,1) & $\phi^{1,0}_{1}(t)\wedge\phi^{0,1}_{2}(t)$, $\phi^{1,0}_{1}(t)\wedge\phi^{0,1}_{3}(t)$, $\phi^{1,0}_{2}(t)\wedge\phi^{0,1}_{1}(t)$\\
 & \\
 &  $\phi^{1,0}_{3}(t)\wedge\phi^{0,1}_{1}(t)$, $\phi^{1,0}_{3}(t)\wedge\phi^{0,1}_{2}(t)$, $\phi^{1,0}_{3}(t)\wedge\phi^{0,1}_{3}(t)$\\
 & \\
(2,1) &  $\phi^{1,0}_{12}(t)\wedge\phi^{0,1}_{1}(t)$, $\phi^{1,0}_{12}(t)\wedge\phi^{0,1}_{2}(t)$, $\phi^{1,0}_{12}(t)\wedge\phi^{0,1}_{3}(t)$\\
 & $\phi^{1,0}_{13}(t)\wedge\phi^{0,1}_{2}(t)$, $\phi^{1,0}_{13}(t)\wedge\phi^{0,1}_{3}(t)$, $\phi^{1,0}_{23}(t)\wedge\phi^{0,1}_{3}(t)$\\
 & \\
(3,1) & $\phi^{1,0}_{123}(t)\wedge\phi^{0,1}_{1}(t)$, $\phi^{1,0}_{123}(t)\wedge\phi^{0,1}_{2}(t)$, $\phi^{1,0}_{123}(t)\wedge\phi^{0,1}_{3}(t)$\\
 & \\
(0,2) & $\phi^{0,1}_{12}(t)$, $\phi^{0,1}_{13}(t)$, $\phi^{0,1}_{23}(t)$\\
 & \\
(1,2) & $\phi^{1,0}_{1}(t)\wedge\phi^{0,1}_{23}(t)$, $\phi^{1,0}_{2}(t)\wedge\phi^{0,1}_{12}(t)$, $\phi^{1,0}_{2}(t)\wedge\phi^{0,1}_{13}(t)$\\
 & $\phi^{1,0}_{3}(t)\wedge\phi^{0,1}_{12}(t)$, $\phi^{1,0}_{3}(t)\wedge\phi^{0,1}_{13}(t)$, $\phi^{1,0}_{3}(t)\wedge\phi^{0,1}_{23}(t)$\\
 & \\
(2,2) & $\phi^{1,0}_{12}(t)\wedge\phi^{0,1}_{12}(t)$, $\phi^{1,0}_{12}(t)\wedge\phi^{0,1}_{13}(t)$, $\phi^{1,0}_{12}(t)\wedge\phi^{0,1}_{23}(t)$\\
 & $\phi^{1,0}_{13}(t)\wedge\phi^{0,1}_{23}(t)$, $\phi^{1,0}_{23}(t)\wedge\phi^{0,1}_{12}(t)$, $\phi^{1,0}_{23}(t)\wedge\phi^{0,1}_{13}(t)$\\
 & \\
(3,2) & $\phi^{1,0}_{123}(t)\wedge\phi^{0,1}_{12}(t)$, $\phi^{1,0}_{123}(t)\wedge\phi^{0,1}_{13}(t)$, $\phi^{1,0}_{123}(t)\wedge\phi^{0,1}_{23}(t)$\\
 & \\
(0,3) & $\phi^{0,1}_{123}(t)$\\
 & \\
(1,3) &  $\phi^{1,0}_{2}(t)\wedge\phi^{0,1}_{123}(t)$, $\phi^{1,0}_{3}(t)\wedge\phi^{0,1}_{123}(t)$\\
 & \\
(2,3) &  $\phi^{1,0}_{12}(t)\wedge\phi^{0,1}_{123}(t)$, $\phi^{1,0}_{23}(t)\wedge\phi^{0,1}_{123}(t)$\\
 & \\
(3,3) & $\phi^{1,0}_{123}(t)\wedge\phi^{0,1}_{123}(t)$\\
\hline
\end{tabular}
\caption{$\Delta''$-Harmonic forms for $X_t$ with $t\neq0$.}\label{tab3}
\end{table}

\begin{table}
\begin{tabular}{l|c}
\hline
$(p,q)$ & Basis for $E^{p,q}_2(X_t)$\\
\hline
 & \\
$(0,0)$ & $1$\\
 & \\
$(1,0)$ & $\phi^{1,0}_1(t)$ \\
 & \\
$(0,1)$ & $\phi^{0,1}_1(t)$\\
 & \\
$(2,0)$ & $\phi^{1,0}_{12}(t)$\\
 & \\
$(1,1)$ & $\phi^{1,0}_{3}(t)\wedge\phi^{0,1}_{1}(t)$, $\phi^{1,0}_{3}(t)\wedge\phi^{0,1}_{2}(t)$\\
 & \\
$(0,2)$ & $\phi^{0,1}_{12}(t)$, $\phi^{0,1}_{13}(t)$, $\phi^{0,1}_{23}(t)$\\
 & \\
$(3,0)$ & $\phi^{1,0}_{123}(t)$\\
 & \\
$(2,1)$ & $\phi^{1,0}_{12}(t)\wedge \phi^{0,1}_{2}(t)$, $\phi^{1,0}_{12}(t)\wedge\phi^{0,1}_{3}(t)$, $\phi^{1,0}_{13}(t)\wedge\phi^{0,1}_{2}(t)$\\
 & \\
$(1,2)$ & $\phi^{1,0}_{1}(t)\wedge\phi^{0,1}_{23}(t)$, $\phi^{1,0}_{2}(t)\wedge\phi^{0,1}_{13}(t)$, $\phi^{1,0}_{3}(t)\wedge\phi^{0,1}_{12}(t)$\\
 & \\
$(0,3)$ & $\phi^{0,1}_{123}(t)$\\
 & \\
$(3,1)$ & $\phi^{1,0}_{123}(t)\wedge\phi^{0,1}_{1}(t)$, $\phi^{1,0}_{123}(t)\wedge \phi^{0,1}_{2}(t)$, $\phi^{1,0}_{123}(t)\wedge \phi^{0,1}_{3}(t)$\\
 & \\
$(2,2)$ & $\phi^{1,0}_{12}(t)\wedge\phi^{0,1}_{13}(t)$, $\phi^{1,0}_{12}(t)\wedge \phi^{0,1}_{23}(t)$\\
 & \\
$(1,3)$ & $\phi^{1,0}_{3}(t)\wedge\phi^{0,1}_{123}(t)$\\
 & \\
$(3,2)$ & $\phi^{1,0}_{123}(t)\wedge\phi^{0,1}_{23}(t)$\\
 & \\
$(2,3)$ & $\phi^{1,0}_{23}(t)\wedge\phi^{0,1}_{123}(t)$\\
 & \\
$(3,3)$ & $\phi^{1,0}_{123}(t)\wedge\phi^{0,1}_{123}(t)$\\
\hline
\end{tabular}
\caption{Basis for $E^{p,q}_2(X_t)$ with $t\neq 0$.}\label{tab4}
\end{table}

\begin{table}
\begin{tabular}{l|c|cc|cc}
\hline
 & & $t=0$ & & $t\neq 0$ & \\
Bi-degree & $b_k$ & $h^{p,q}(0)$ & $e^{p,q}_2(0)$ & $h^{p,q}(t)$ & $e^{p,q}_2(t)$\\
\hline
 & & & & & \\
$(0,0)$ & 1 & 1 & 1 & 1 & 1\\ 
\hline
 & & & & & \\ 
$(1,0)$ & & 3 & 1 & 2 & 1\\
 & 2 & & & & \\ 
$(0,1)$ & & 3 & 1 & 3 & 1\\
\hline
 & & & & & \\ 
$(2,0)$ & & 3 & 1 & 2 & 1 \\
 & & & & & \\ 
$(1,1)$ & 5 & 9 & 3 & 6 & 2\\
 & & & & & \\ 
$(0,2)$ & & 3 & 1 & 3 & 3\\
\hline
 & & & & & \\ 
$(3,0)$ & & 1 & 1 & 1 & 1\\
 & & & & & \\ 
$(2,1)$ & & 9 & 3 & 6 & 3\\
 & 8 & & & & \\ 
$(1,2)$ & & 9 & 3 & 6 & 3\\
 & & & & & \\ 
$(0,3)$ & & 1 & 1 & 1 & 1\\
\hline
 & & & & & \\ 
$(3,1)$ & & 3 & 1 & 3 & 3\\
 & & & & & \\ 
$(2,2)$ & 5 & 9 & 3 & 6 & 2\\
 & & & & & \\ 
$(1,3)$ & & 3 & 1 & 2 & 1\\
\hline
 & & & & & \\ 
$(3,2)$ & & 3 & 1 & 3 & 1\\
 & 2 & & & & \\ 
$(2,3)$ & & 3 & 1 & 2 & 1\\
\hline
 & & & & & \\ 
$(3,3)$ & 1 & 1 & 1 & 1 & 1\\
\hline 
\end{tabular}
\caption{Comparison between dimensions of $E^{\bullet,\bullet}_r(X_0)$ and $E^{\bullet,\bullet}_r(X_t)$, with $r=1,2$ and $t\neq 0$.}\label{tab5}
\end{table}

\end{document}